# Technical details regarding infinite families of crossing-critical graphs with prescribed average degree and crossing number

Drago Bokal (drago.bokal@imfm.uni-lj.si)

Institute of mathemathics, physics and mechanics
Department of Mathematics
Jadranska 19
Ljubljana, Slovenia

## Abstract

This *Mathematica* notebook is presented as a supplement of the paper [1]. The main result is the following theorem:

**For every interval $I = [r_1, r_2]$, $r_1, r_2 \in (3, 6)$ there exists an integer $N_I > 0$, such that for every rational number $r \in I$, and for every integer $k > N_I$, there exists an infinite family of crossing critical graphs, all having average degree $r$ and crossing number $k$.**

The statement holds for $N_I = $ `Ceiling`[`Max`[$f(r_1), f(r_2)$]], where $f(x) = 240 + \frac{512}{(6-x)^2} + \frac{224}{6-x} + \frac{25}{16(x-3)^2} + \frac{40}{x-3}$.

This notebook contains technical details omitted in the paper and can be used as a hint of how to rigorously verify the constraints that are imposed on the parameters in the main construction of the paper. The reader may either use *Mathematica* or some other software to verify the listed claims, or may derive them in a more clear, oldfashioned way. As the notebook is not self-contained, it is advisable to read the paper before.

## 1. Some commands

The command `Off` disables the warning that repeatedly pops up due to similar names of the symbols used. `FS` is a shortening of the command `FullSimplify`, either without or with assumptions to be respected during simplification. The command `AverageDegree` is used to compute the average degree, given the partial degree sequence of the graph. `DegreeSequence3456` encapsulates the quadruple containing the number of vertices with degrees 3, 4, 5, and 6, respectively. Given two such partial degree sequences, the command `Zip3` computes the partial degree sequence of the graph, obtained as the zip product at vertices of degree 3 of two graphs with respective degree sequences.

```
In[1]:= Off[General::spell1];
       FS[e_] := FullSimplify[e];
       FS[e_, a_] := FullSimplify[e, a];
       AverageDegree[d_DegreeSequence3456] :=
         FS[(3 d[[1]] + 4 d[[2]] + 5 d[[3]] + 6 d[[4]]) / (d[[1]] + d[[2]] + d[[3]] + d[[4]])];
       Zip3[d1_DegreeSequence3456, d2_DegreeSequence3456] := DegreeSequence3456[
          d1[[1]] + d2[[1]] - 2, d1[[2]] + d2[[2]], d1[[3]] + d2[[3]], d1[[4]] + d2[[4]]];
```

# 2. Graphs with average degree close to 3

The following section presents the construction of the crossing number critical graphs with average degree greater, but arbitratily close to 3. The proof of the crossing number relies on the presence of the twisted pueblo substructure in the graph and is given in the aforementioned paper. A graph in the family $S(n, m, c)$ is obtained from the sequence of tiles $S_n, S_n', S_n, S_n',..., S_n', S_n'', S_n, S_n', S_n, ...,S_n'$ of odd length $m$ strictly more than $4\binom{n}{2}$, where $S_n$ is the staircase tile (Figure A), $S_n'$ is the inverted staircase tile (Figure B) and $S_n''$ is the twisted staircase tile (Figure C). Precisely $c$ thick edges (see the next Figure) are contracted after joining the tiles. For each of the different ways in which the contraction can be done, a different graph of the family $S(n, m, c)$ is obtained.

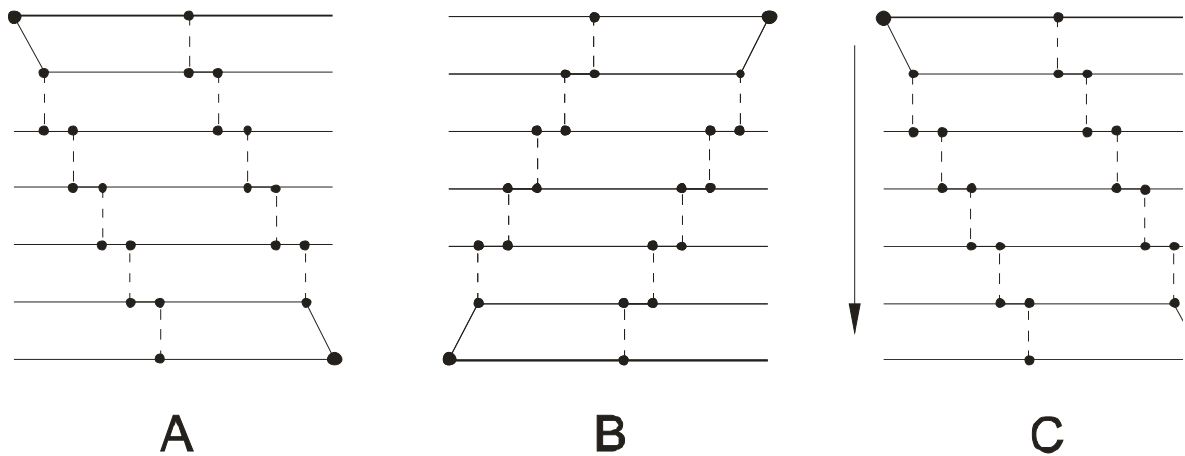

A          B          C

The following figure illustrates that the tile crossing number of a join of a twisted staircase tile is at most $\binom{n}{2} - 1$. The drawing can be augmented so that any of the dashed edges is crossed, thus removing any edge from the tile decreases its crossing number. Note that contracting the thick edges does not affect the crossing number.

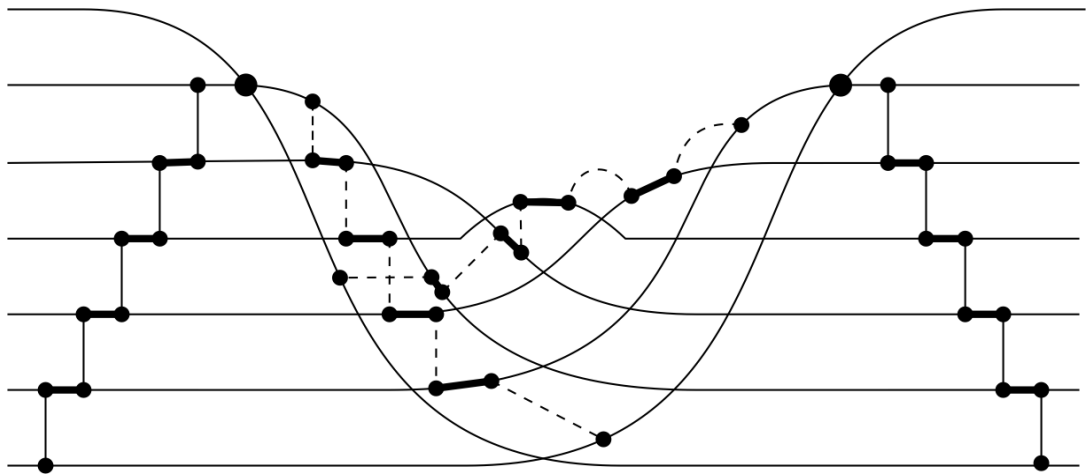

The graphs from $S(n, m, c)$ have $4m(n-2) - 2c$ vertices of degree 3 and $m + c$ vertices of degree 4. Each of them has $2m(n-3)$ thick edges that can be contracted. This number is computed by the command `Thick`. The command `GraphS` produces the 3, 4, 5, 6 partial degree sequence, which is used to compute the average degree. The formulas account for the fact that $m$ is odd.

```
In[6]:=  Thick[n_, m_] := 2 m (n - 3);
         GraphS[n_, m_, c_] := DegreeSequence3456[4 m (n - 2) - 2 c, m + c, 0, 0];
         m = 2 mp + 1; (* m has to be odd. *)
         AverageDegree[GraphS[n, m, c]] // FS
```

$$Out[9]= \quad 3 + \frac{1 + c + 2\, mp}{-c + (1 + 2\, mp)\,(-7 + 4\, n)}$$

By choosing $n$ sufficiently large, the average degree of the graphs in $S(n, m, 0)$ can be arbitrarily close to 3. Due to the presence of the twisted pueblo substructure in the graph $S(n, m, c)$, for $m$ larger than $4\binom{n}{2}$, the crossing number of the graphs in $S(n, m, c)$ is $\binom{n}{2} - 1$.

```
In[10]:= CrnS[n_] := Binomial[n, 2] - 1;
         CrnS[n]
```

$$Out[11]= \quad -1 + \frac{1}{2}(-1 + n)\,n$$

The following constraints on the parameters of the family $S(n, m, c)$ have to be satisfied:

```
In[12]:= ConstraintsS[n_, m_, c_] := {n ∈ Integers, n ≥ 3,
            (m - 1) / 2 ∈ Integers, m > 4 CrnS[n], c ∈ Integers, c ≥ 0, c ≤ Thick[n, m]};

         ConstraintsS[n, m, c]
```

$$Out[13]= \quad \Big\{n \in \text{Integers}, n \geq 3,\ mp \in \text{Integers},$$
$$1 + 2\,mp > 4\left(-1 + \frac{1}{2}(-1 + n)\,n\right),\ c \in \text{Integers}, c \geq 0, c \leq 2(1 + 2\,mp)(-3 + n)\Big\}$$

## Proof of proposition 15

If $0 < a < b$ and $a + b$ is odd, one may construct an infinite family of graphs with average degree $3 + \frac{a}{b}$ just by carefully choosing the parameters of the family $S(n, m, c)$. Choose $n \geq \max\left(\frac{5b-a}{2(b-a)}, \frac{7a+b}{4a}, 3\right)$, $k > n^2$ and set
$m(t) = (2t + 1)(a + b)$
$c(t) = (2t + 1)((4n - 7)a - b)$.

```
In[14]:= AssumptionsS = {n ≥ 4, n ≥ (5 b - a)/(2 (b - a)), n ≥ (7 a + b)/(4 a), t ≥ n^2,
           b > a, a > 0, a ∈ Integers, b ∈ Integers, n ∈ Integers, t ∈ Integers}
         MM[t_] := (2 t + 1) (a + b);
         CC[t_] := (2 t + 1) ((4 n - 7) a - b);
```

$$Out[14]= \quad \Big\{n \geq 4,\ n \geq \frac{-a + 5b}{2(-a + b)},\ n \geq \frac{7a + b}{4a},\ t \geq n^2,\ b > a,$$
$$a > 0,\ a \in \text{Integers},\ b \in \text{Integers},\ n \in \text{Integers},\ t \in \text{Integers}\Big\}$$

Then the graphs in $S(n, m(t), c(t))$ are critical graphs with crossing number $\binom{n}{2} - 1$. We need to verify the following constraints on the values of the parameters:

The value of *c* must be nonnegative:

*In[17]:=* `FullSimplify[CC[t] ≥ 0, AssumptionsS]`

*Out[17]=* True

The value of *c* must be smaller or equal to the overall number of thick edges:

*In[18]:=* `FullSimplify[CC[t] ≤ Thick[n, MM[t]], AssumptionsS]`

*Out[18]=* True

The value of *m* must be greater than $4 \operatorname{cr}(S(n, m, c))$:

*In[19]:=* `FullSimplify[MM[t] > 4 CrnS[n], AssumptionsS]`

*Out[19]=* True

As demonstrated, all the constraints for S(n,m(t),c(t)) to be critical graphs are satisfied. The only property needed to verify is the desired average degree:

*In[20]:=* `AverageDegree[GraphS[n, MM[t], CC[t]]]`

*Out[20]=* $3 + \dfrac{a}{b}$

# Graphs with average degree close to 6

Construction of the crossing number critical graphs $H(w, s)$ with average degree smaller, but arbitratily close to 6, is presented in this section. As opposed to the previous construction, there is only one graph $H(w, s)$ for given values of the parameters. The graph $H(w, s)$ is obtained by joining a sequence of tiles $H_w, H_w, \ldots, H_w, H_w', H_w$ of length $s$, where the tile $H_w$, for $w = 1$, is presented in the following Figure, and $H_w'$ is the twist of this tile. The letters in the Figure denote traversing paths of the tile, which are used to form twisted pairs in the twisted tile. The number of such twisted pairs determines a lower bound on the tile crossing number of the join of the sequence of tiles. As each path extends when the tiles are joined, the lower bound on the tile crossing number holds for joins of the sequences of any length.

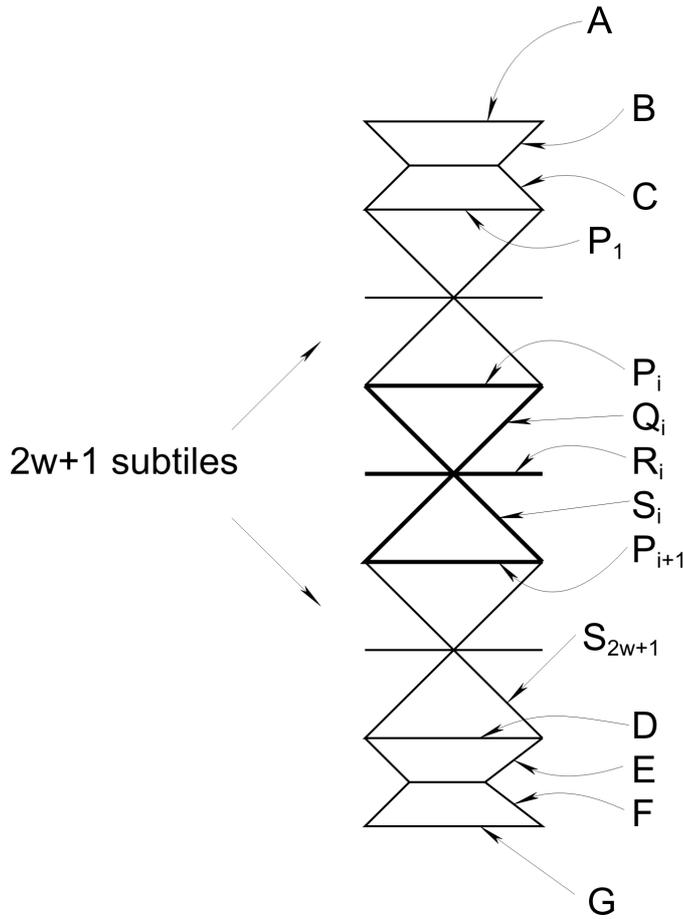

Each tile of the graph contributes 4 vertices of degree 3, 2 vertices of degree 4 and $4w + 3$ vertices of degree 6. The 3, 4, 5, 6 degree sequence of the graph is computed using the command GraphH.

```
In[21]:=  GraphH[w_, s_] := DegreeSequence3456[4 s, 2 s, 0, (4 w + 3) s];
          AverageDegree[GraphH[w, s]] // FS
```

$$Out[22]= \ 6 - \frac{16}{9 + 4w}$$

For $w$ sufficiently large, the average degree is arbitarily close to 6. The crossing number of the graphs $H(w, s)$ is computed in what follows. First, the lower bound is established based on the ideas od Pinontoan and Richter [3]. For large $s$, the lower bound is obtained by counting pairs of traversing paths $(P, Q)$, such that:

- in any pair $(P, Q)$ the paths $P$ and $Q$ are vertex disjoint and twisted in the twist of the tile presented above,
- for any two distinct pairs $(P, Q)$ and $(P', Q')$, either $P$ and $P'$ or $Q$ and $Q'$ are edge disjoint.

A family with these properties is described in terms of pairs of paths we denoted with the letters $A$, $B$, $C$, $D$, $P_i$, $Q_i$, $R_i$, $S_i$ in the above Figure. The pairs $(P, Q)$ are considered lexicographically by the letters indicating the paths in them. A pair $(P, Q)$ is called *valid*, if the set of all the pairs, mentioned before, preserves the above two properties after adding the pair $(P, Q)$ to it.

The path $A$ forms a valid pair with each of $C$, $D$, $E$, $G$, and $P_i$, $Q_i$, $R_i$, $S_i$, for all $i$, thus the number of pairs is:

```
In[23]:=  PairsA[w_] := 4 + 4 (2 w + 1) // FS
          PairsA[w]
```

$$Out[24]= \ 8 (1 + w)$$

The path $B$ forms a disjoint pair with $D$, $E$, $G$, and $P_i$, $Q_i$, $R_i$, $S_i$, for all $i$, thus the number of pairs is:

```
In[25]:= PairsB[w_] := 3 + 4 (2 w + 1) // FS
        PairsB[w]
```

*Out[26]=*   7 + 8 w

$C$ has an edge in common with $B$, thus the only possible pair is $(A, C)$, which is already counted in `PairsA`. $D$ forms a disjoint pair with the paths $F$, $G$, the paths $P_i$, $Q_i$, $R_i$, for all $i$, and the paths $S_i$, for all but one $i$. Thus,

```
In[27]:= PairsC[w_] := 0;
        PairsD[w_] := 2 + 4 (2 w + 1) - 1 // FS
        PairsC[w]
        PairsD[w]
```

*Out[29]=*   0

*Out[30]=*   5 + 8 w

The path $E$ forms a disjoint pair with $G$, the paths $P_i$, $Q_i$, $R_i$, for all $i$, and the paths $S_i$, for all but one $i$. Thus,

```
In[31]:= PairsE[w_] := 1 + 4 (2 w + 1) - 1 // FS
        PairsE[w]
```

*Out[32]=*   4 + 8 w

The path $F$ forms a valid disjoint pair only with $S_{2w+1}$.

```
In[33]:= PairsF[w_] := 1;
        PairsF[w]
```

*Out[34]=*   1

The path $G$ forms valid disjoint pairs with the paths $P_i$, $Q_i$, $R_i$, and $S_i$ for all $i$. Thus,

```
In[35]:= PairsG[w_] := 4 (2 w + 1)
        PairsG[w]
```

*Out[36]=*   4 (1 + 2 w)

For a given $i$, the path $P_i$ forms a valid disjoint pair with $R_j$ and $S_j$ for all $j \geq i$ and with $P_j$, $Q_j$ for all $j > i$. Thus,

```
In[37]:= PairsPi[i_] := 2 (2 w + 2 - i) + 2 (2 w + 1 - i) // FS;
        PairsP[w_] := Sum[PairsPi[i], {i, 1, 2 w + 1}] // FS
        PairsP[w]
```

*Out[39]=*   $2 (1 + 2 w)^2$

For a given $i$, the paths $Q_i$ and $R_i$ form a valid disjoint pair $P_j$, $Q_j$, $R_j$, $S_j$ for all $j > i$. Thus,

```
In[40]:= PairsQi[i_] := 4 (2 w + 1 - i) // FS;
        PairsQ[w_] := Sum[PairsQi[i], {i, 1, 2 w + 1}] // FS
        PairsR[w_] := PairsQ[w];
        PairsR[w]
```

*Out[43]=*   4 w (1 + 2 w)

For a given $i \leq 2w$, the path $S_i$ forms a valid pair with $R_j$ and $S_j$ for $j > i$ and with $P_j$, $Q_j$ for $j > i + 1$. The path $S_{2w+j}$ forms no valid pair.

```
In[44]:=  PairsSi[i_] := 2 (2 w + 1 - i) + 2 (2 w - i);
          PairsS[w_] := Sum[PairsSi[i], {i, 1, 2 w}]
          PairsS[w]
```

*Out[46]=*  $8 w^2$

The tile crossing number of the join of the sequence of tiles is bounded from below by the overall number of pairs, which is:

```
In[47]:=  CrnH[w_] := PairsA[w] + PairsB[w] + PairsC[w] + PairsD[w] + PairsE[w] +
             PairsF[w] + PairsG[w] + PairsP[w] + PairsQ[w] + PairsR[w] + PairsS[w] // FS;

          CrnH[w]
```

*Out[48]=*  $31 + 8 w (7 + 4 w)$

The following Figure shows a tile drawing of the graph for $w = 0$ with 31 crossings, denoted as hollow circles. The drawing can be easily generalized to a drawing for $w > 0$ with $32 w^2 + 56 w + 31$ crossings. It is clear that whichever edge is removed from the tile, the crossing number is decreased.

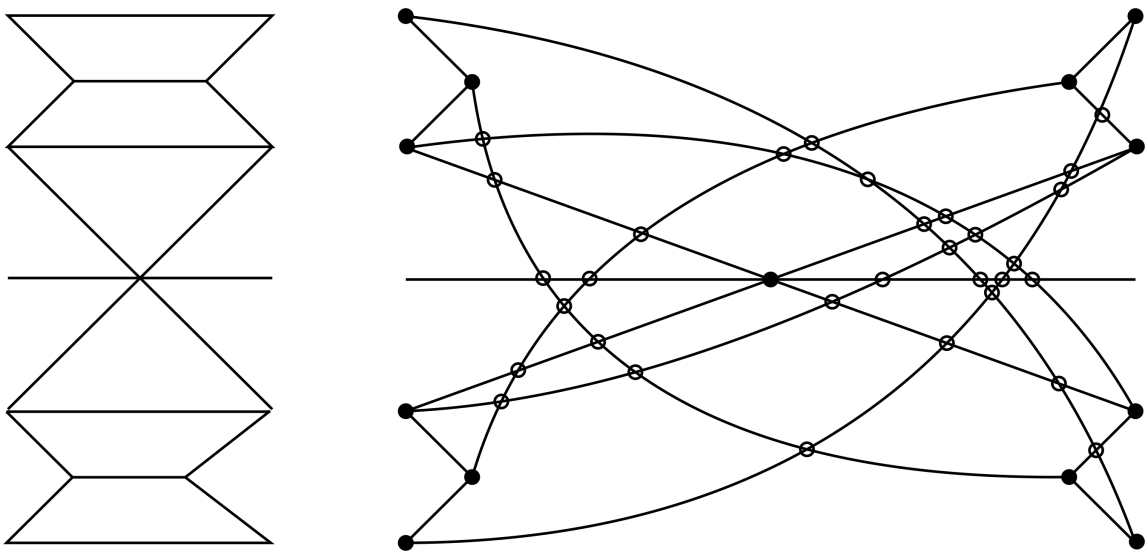

By construction, the graph $H(w, s)$ will have the stated crossing number for $s$ larger than $4 (32 w^2 + 56 w + 32)$. The following lists the constraints on parameters $w$ and $s$.

```
In[49]:=  ConstraintsH[w_, s_] = {w ∈ Integers, w ≥ 0, s ∈ Integers, s > 4 CrnH[w]};
          ConstraintsH[w, s]
```

*Out[50]=*  {w ∈ Integers, w ≥ 0, s ∈ Integers, s > 4 (31 + 8 w (7 + 4 w))}

# Adapting graphs

Adapting graphs $R(p, q)$ are the graphs used to join a graph with average degree close to 6 and a graph with average degree close to 3, thus yielding a graph with desired average degree between 3 and 6. The three graphs are joined using zip product at vertices of degree 3. The zip product compromises the pattern in the graphs $H$ and $S$, and the role of the graph $R$ is to compensate for this breach. The graph $R$ is also used to obtain the desired crossing number of the resulting graph.

The family of graphs $R(p, q)$ is obtained inductively in the following way: $R(1, 0)$ contains only the complete bipartite graph $K_{3,3}$. Graphs in the family $R(p, 0)$ are obtained by the zip product of some graph $R$ from $R(p-1, 0)$ and $K_{3,3}$. The graphs of $R(p, q)$ are obtained as a zip product of a graph $R$ from $R(p, q-1)$ and the graph $K_{3,5}$ at a vertex of degree 3.

The graph in $R(1, 0)$ has 6 vertices of degree 3. A zip product with $K_{3,3}$ adds 4 vertices of degree 3, thus $R(p, 0)$ has $4p + 2$ vertices of degree 3. A zip product with $K_{3,5}$ at a vertex of degree 3 adds 3 vertices of degree 3 and 3 vertices of degree 5. Thus the graph $R(p, q)$ has $4p + 3q + 2$ vertices of degree 3, and $3q$ vertices of degree 5. The command GraphR produces the 3, 4, 5, 6 degree sequence of the graph $R(p, q)$.

*In[51]:=* **GraphR[p_, q_] = DegreeSequence3456[4 p + 3 q + 2, 0, 3 q, 0];**
**AverageDegree[GraphR[p, q]] // FS**

*Out[52]=* $\dfrac{3 + 6p + 12q}{1 + 2p + 3q}$

Average degree of the graphs in $R(p, q)$ is between 3 (for high $p$) and 4 (for high $q$):

*In[53]:=* **AverageDegree[GraphR[p, q]] - 3 // FS**
**4 - AverageDegree[GraphR[p, q]] // FS**

*Out[53]=* $\dfrac{3q}{1 + 2p + 3q}$

*Out[54]=* $\dfrac{1 + 2p}{1 + 2p + 3q}$

Crossing number of the graphs $R(p, q)$ is computed using the known crossing numbers of $K_{3,3}$ (being 1) and the known crossing number of $K_{3,5}$ (being 4, due to Kleitman [2]). As the zip product preserves the crossing number, the crossing number of $R(p, q)$ is

*In[55]:=* **CrnR[p_, q_] := p + 4 q;**
**CrnR[p, q]**

*Out[56]=* p + 4 q

The following constraints are imposed on the parameters of the family $R(p, q)$:

*In[57]:=* **ConstraintsR[p_, q_] = {p ∈ Integers, p ≥ 1, q ∈ Integers, q ≥ 1};**
**ConstraintsR[p, q]**

*Out[58]=* {p ∈ Integers, p ≥ 1, q ∈ Integers, q ≥ 1}

# Combining the graphs

The graphs in the family $\Gamma(n, m, c, w, s, p, q)$ are obtained as a zip product of some graph $S \in S(n, m, c)$, $H = H(w, s)$, and some $R \in R(p, q)$. Both zip products are performed at vertices of degree 3. The choice of vertices and the order in which this is done are irrelevant; the resulting graph always has the same average degree and crossing number. The degree sequence, the average degree, the crossing number, and the constraints on the parameters of the family $\Gamma(n, m, c, w, s, p, q)$ are as follows:

```
In[59]:= GraphG[n_, m_, c_, w_, s_, p_, q_] :=
           Zip3[Zip3[GraphS[n, m, c], GraphH[w, s]], GraphR[p, q]];
         CrnG[n_, m_, c_, w_, s_, p_, q_] := CrnS[n] + CrnH[w] + CrnR[p, q];
         ConstraintsG[n_, m_, c_, w_, s_, p_, q_] =
           Union[ConstraintsS[n, m, c], ConstraintsH[w, s], ConstraintsR[p, q]];
         AverageDegree[GraphG[n, m, c, w, s, p, q]] // FS
         CrnG[n, m, c, w, s, p, q]
         ConstraintsG[n, m, c, w, s, p, q]
```

$$Out[62]= 6 + \frac{4\,(7 + c + mp\,(11 - 6\,n) - 3\,n - 3\,p - 3\,q - 4\,s)}{-9 - c + 4\,n + 2\,mp\,(-7 + 4\,n) + 4\,p + 6\,q + 9\,s + 4\,s\,w}$$

$$Out[63]= 30 + \frac{1}{2}\,(-1 + n)\,n + p + 4\,q + 8\,w\,(7 + 4\,w)$$

$$Out[64]= \{c \in \text{Integers}, mp \in \text{Integers}, n \in \text{Integers}, p \in \text{Integers},$$
$$q \in \text{Integers}, s \in \text{Integers}, w \in \text{Integers}, 1 + 2\,mp > 4\,\left(-1 + \frac{1}{2}\,(-1 + n)\,n\right),$$
$$s > 4\,(31 + 8\,w\,(7 + 4\,w)), c \geq 0, n \geq 3, p \geq 1, q \geq 1, w \geq 0, c \leq 2\,(1 + 2\,mp)\,(-3 + n)\}$$

Let $r = 3 + \frac{a}{b}$, $3 < r < 6$ be given. As the parameters used in construction need to be integers, several integer divisions are performed to obtain auxiliary values used in defining the final parameters. Let $b = b'\,a + b_r$, $b' = 4\,b'' + b_r'$, and $4\,b = \bar{b}(3\,b - a) + \bar{b}_r$. In these equalities, $b_r$, $b_r'$ and $\bar{b}_r$ are the reminders of the integer divisions and are thus nonnegative integers, smaller than the respective divisors. As the parameters are nonnegative, also the coefficients are nonnegative integers.

Now let $N = 8\,\bar{b}\,(4\,\bar{b} + 7) + 5\,(b'' + 4)\,(5\,b'' + 12)$ and let $k > N$. Perform another division and let $k'\,(2\,b'' + 5) + k_r = k - \frac{1}{2}\,b''(b'' + 5) - 8\,\bar{b}\,(4\,\bar{b} + 7) - 34$.

The following choice of the parameters of $\Gamma(n, m, c, w, s, p, q)$ assures that every graph $G \in \Gamma(n, m, c, w, s, p, q)$ is a $k$-critical graph with crossing number $k$ and average degree $r$. First we choose an arbitrary integer $t > k$. Then we set
$n = 4 + b''$,
$m = 2\,t\,(27\,b - 9\,a - 4\,\bar{b}_r) - 2\,k' + 3$,
$c = 2\,k' - 12\,b'' - 6\,k_r - 33$,
$w = \bar{b}$,
$s = 2\,t\,(a\,(4\,b'' + 9) - b)$,
$p = k - (\frac{1}{2}\,b''(b'' + 23) + 8\,\bar{b}\,(4\,\bar{b} + 7) + 4\,k_r + 56)$,
$q = 2\,b'' + k_r + 5$.

According to the first three integer divisions and the lower bound on $k$ and $t$, the following assumptions on the values of the parameters hold:

```
In[65]:= b == a bp + br;
        bp == 4 bpp + bpr;
        4 b == bb (3 b - a) + bbr;
        NN = 8 bb (4 bb + 7) + 5 (bpp + 4) (5 bpp + 12);
        AssumptionsG =
          {a ∈ Integers, b ∈ Integers, a ≥ 1, 3 b > a, br ≥ 0, br ∈ Integers, bp ∈ Integers, bp ≥ 0,
           br < a, br ≤ b, bpp ∈ Integers, bpr ∈ Integers, bpp ≥ 0, bpr ≥ 0, bpr ≤ 3, bpr ≤ bp,
           bb ∈ Integers, bbr ∈ Integers, bb ≥ 0, bbr ≥ 0, bbr < 3 b - a, bbr ≤ 4 b, k > NN, t > k};
```

The integer $k - \frac{1}{2} b''(b'' + 5) - 8 \bar{b}(4\bar{b} + 7) - 34$ participating in the fourth division is positive:

```
In[70]:= FS[k - bpp (bpp + 5) / 2 - 8 bb (4 bb + 7) - 34 > 0, AssumptionsG]

Out[70]= True
```

The fourth integer division allows the following new assumptions:

```
In[71]:= k - bpp (bpp + 5) / 2 - 8 bb (4 bb + 7) - 34 == kp (2 bpp + 5) + kr;
        AssumptionsG = Union[AssumptionsG, {kr ∈ Integers, kp ∈ Integers, kp ≥ 0,
            kr ≥ 0, kr < (2 bpp + 5), kr ≤ k - bpp (bpp + 5) / 2 - 8 bb (4 bb + 7) - 34}];
```

The constraints of the construction are verified in the sequel.

```
In[73]:= ConstraintsS[n, m, c]

Out[73]= {n ∈ Integers, n ≥ 3, mp ∈ Integers,
         1 + 2 mp > 4 (-1 + 1/2 (-1 + n) n), c ∈ Integers, c ≥ 0, c ≤ 2 (1 + 2 mp) (-3 + n)}
```

The value of parameter $n = 4 + b''$ is clearly an integer. The constraint $n ≥ 3$ is satisfied:

```
In[74]:= n = bpp + 4;
        FS[n ≥ 3, AssumptionsG]

Out[75]= True
```

The parameter $m = 2t(27 b - 9a - \bar{b}_r) - 2k' + 3$ is an odd number, as desired. The following establishes that it respects the imposed lower bound:

```
In[76]:= m = 2 t (27 b - 9 a - 4 bbr) - 2 kp + 3;
        FS[m > 4 (-1 + 1/2 (-1 + n) n), AssumptionsG]

Out[77]= 54 b t > 17 + 2 bpp (7 + bpp) + 2 kp + 18 a t + 8 bbr t
```

The algorithm needs some hints to verify the constraint. In the following, both sides of the inequality are put together and the coefficient at *t* is determined:

```
In[78]:= 54 b t - (27 + 3 bpp (7 + bpp) + 2 kp + 18 a t + 8 bbr t) // ExpandAll // Collect[#, t] &

Out[78]= -27 - 21 bpp - 3 bpp² - 2 kp + (-18 a + 54 b - 8 bbr) t
```

The following verifies that the coefficient at *t* is at least 1.

```
In[79]:= FS[(-18 a + 54 b - 8 bbr) > 1, AssumptionsG]

Out[79]= True
```

As $t > k$, the following establishes the validity of the lower bound on $m$. In order that the algorithm is able to simplify the expression, the value of $k'$ needs to be replaced by the formula obtained from the respective integer division.

*In[80]:=* `Solve[k - bpp (bpp + 5) / 2 - 8 bb (4 bb + 7) - 34 == kp (2 bpp + 5) + kr, kp]`

*Out[80]=* $\left\{\left\{kp \to \dfrac{-68 - 112\,bb - 64\,bb^2 - 5\,bpp - bpp^2 + 2\,k - 2\,kr}{2\,(5 + 2\,bpp)}\right\}\right\}$

*In[81]:=* `FS[-27 - 21 bpp - 3 bpp^2 - 2 kp + k > 0 /.`
`        kp ->` $\dfrac{-68 - 112\,bb - 64\,bb^2 - 5\,bpp - bpp^2 + 2\,k - 2\,kr}{2\,(5 + 2\,bpp)}$ `, AssumptionsG]`

*Out[81]=* `True`

The value of parameter $c = 2\,k' - 12\,b'' - 6\,k_r - 33$ is clearly an integer. It has to respect a lower bound $c \geq 0$ and an upper bound $c \leq 2\,m\,(n - 3)$.

*In[82]:=* `c = 2 kp - 12 bpp - 6 kr - 33;`
`        FS[c >= 0 /. kp ->` $\dfrac{-68 - 112\,bb - 64\,bb^2 - 5\,bpp - bpp^2 + 2\,k - 2\,kr}{2\,(5 + 2\,bpp)}$ `, AssumptionsG]`

*Out[83]=* `True`

The upper bound is slightly more demanding:

*In[84]:=* `FS[c <= 2 m (n - 3), AssumptionsG]`

*Out[84]=* `6 kp + 4 bpp kp + 4 (9 a - 27 b + 4 bbr) (1 + bpp) t <= 39 + 18 bpp + 6 kr`

The coefficient at $t$ is smaller than $-2$:

*In[85]:=* `FS[4 (9 a - 27 b + 4 bbr) (1 + bpp) < -2 /. bb ->` $\dfrac{-4\,b + bbr}{a - 3\,b}$ `, AssumptionsG]`

*Out[85]=* `True`

As $t > k$, the following establishes the upper bound on $c$:

*In[86]:=* `FS[6 kp + 4 bpp kp - 2 k <= 39 + 18 bpp + 6 kr /.`
`        kp ->` $\dfrac{-68 - 112\,bb - 64\,bb^2 - 5\,bpp - bpp^2 + 2\,k - 2\,kr}{2\,(5 + 2\,bpp)}$ `, AssumptionsG]`

*Out[86]=* `True`

Thus, all the constraints that a graph $S \in S(n, m, c)$ is a critical graph, are satisfied. The graph $H = H(w, s)$ is analyzed in the following. Its parameters have to respect the following constraints:

*In[87]:=* `ConstraintsH[w, s]`

*Out[87]=* `{w ∈ Integers, w >= 0, s ∈ Integers, s > 4 (31 + 8 w (7 + 4 w))}`

Recall that $w = \overline{b}$.

*In[88]:=* `w = bb;`
`        FS[w >= 0, AssumptionsG]`

*Out[89]=* `True`

The lower bound on $s = 2t(a(4b''+9) - b)$ is somewhat more difficult:

```
In[90]:= s = 2 t (a (4 bpp + 9) - b);
        FS[s > 4 (31 + 8 w (7 + 4 w)) /. bb → (-4 b + bbr)/(a - 3 b), AssumptionsG]
```

$$Out[91]= -2(b - a(9 + 4\,bpp))\,t > 4\left(31 - \frac{8(4b - bbr)(7a - 37b + 4\,bbr)}{(a - 3b)^2}\right)$$

The coefficient at $t$ is at least 4:

```
In[92]:= FS[-2 (b - a (9 + 4 bpp)) > 4 /. bpp → (b - a bpr - br)/(4 a), AssumptionsG]

Out[92]= True
```

As $t > k$, the following establishes the constraint on $s$:

```
In[93]:= FS[4 k > 4 (31 - (8 (4 b - bbr) (7 a - 37 b + 4 bbr))/(a - 3 b)^2) /. bbr → 4 b + a bb - 3 b bb, AssumptionsG]

Out[93]= True
```

The only assumptions that are left to be verified are those on the graph $R \in R(p, q)$.

```
In[94]:= ConstraintsR[p, q]

Out[94]= {p ∈ Integers, p ≥ 1, q ∈ Integers, q ≥ 1}
```

The constraint on $p = k - (\frac{1}{2} b''(b''+23) + 8\overline{b}(4\overline{b}+7) + 4k_r + 56)$ is established in the following:

```
In[95]:= p = k - (bpp (bpp + 23) / 2 + 8 bb (4 bb + 7) + 4 kr + 56);
        FS[p ≥ 1, AssumptionsG]

Out[96]= True
```

Also the constraint on $q = 2b'' + k_r + 5$ does not present any difficulties:

```
In[97]:= q = 2 bpp + kr + 5;
        FS[q ≥ 1, AssumptionsG]

Out[98]= True
```

The crossing number of the graphs in $\Gamma$ is as desired:

```
In[99]:= FS[CrnG[n, m, c, w, s, p, q]]

Out[99]= k
```

So is the average degree:

*In[100]:=*
```
FS[AverageDegree[GraphG[n, m, c, w, s, p, q]]] /.
  kp → (-68 - 112 bb - 64 bb^2 - 5 bpp - bpp^2 + 2 k - 2 kr)/(2 (5 + 2 bpp)) /.
  bbr → 4 b + a bb - 3 b bb, AssumptionsG]
```

*Out[100]=*
$$3 + \frac{a}{b}$$

A closer observation of the bound *N* reveals that the lower bound on the crossing number can be replaced by a function of $x = 3 + \frac{a}{b}$.

*In[101]:=*
```
NN
```

*Out[101]=*
$$8\,bb\,(7 + 4\,bb) + 5\,(4 + bpp)\,(12 + 5\,bpp)$$

*In[102]:=*
```
b == a bp + br;
bp == 4 bpp + bpr;
Solve[4 b == bb (3 b - a) + bbr, bb]
```

*Out[104]=*
$$\left\{\left\{bb \to \frac{-4\,b + bbr}{a - 3\,b}\right\}\right\}$$

As $\bar{b}_r$ is at least 0 and at most $3b - a$, $\bar{b}$ lies in the interval $(\frac{4}{6-x} - 1, \frac{4}{6-x}]$.

*In[105]:=*
```
Solve[b == a (4 bpp + bpr) + br, bpp]
```

*Out[105]=*
$$\left\{\left\{bpp \to \frac{b - a\,bpr - br}{4\,a}\right\}\right\}$$

Similarly, as $b_r$ is less than *a* and at least 0, and as $b_r'$ is at most 3 and at least 0, *b″* lies in the interval $\left(\frac{1}{4(x-3)} - 1, \frac{1}{4(x-3)}\right]$. This demonstrates that the following bound can replace the bound *N*:

*In[106]:=*
```
NN /. {bpp -> 1/(4 (x - 3)), bb -> 4/(6 - x)} // Apart
```

*Out[106]=*
$$240 + \frac{512}{(-6 + x)^2} - \frac{224}{-6 + x} + \frac{25}{16\,(-3 + x)^2} + \frac{40}{-3 + x}$$

*In[107]:=*
```
Bnd[x_] := 240 + 512/(-6 + x)^2 - 224/(-6 + x) + 25/(16 (-3 + x)^2) + 40/(-3 + x);
```

*In[108]:=*
```
Plot[Bnd[x], {x, 3.01, 5.9}]
```

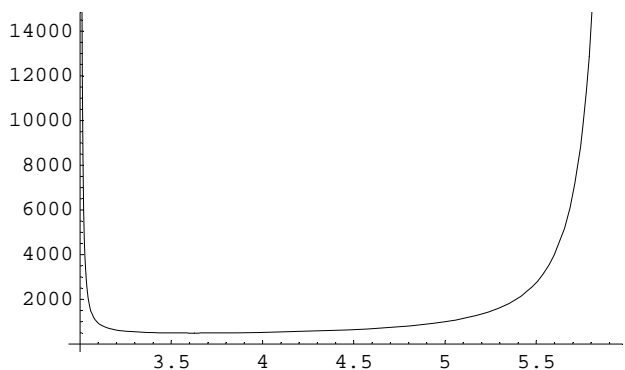

*Out[108]=*
     - Graphics -

*In[109]:=*
```
D[Bnd[x], {x, 2}] // FS
```

*Out[109]=*
$$\frac{3072}{(-6+x)^4} - \frac{448}{(-6+x)^3} + \frac{75}{8\,(-3+x)^4} + \frac{80}{(-3+x)^3}$$

*In[110]:=*
```
FS[D[Bnd[x], {x, 2}] > 0, {3 < x, x < 6}]
```

*Out[110]=*
     True

The function has a pole at $x = 0$ and $x = 3$ and is convex in the interval $(0, 3)$, thus also in any interval $I$ contained in $(0, 3)$. Its supremum in $I$ is attained at one of the boundary points of $I$.

**For every interval $I = [r_1, r_2]$, $r_1, r_2 \in (3, 6)$ there exists an integer $N_I > 0$, such that for every rational number $r \in I$, and for every integer $k > N_I$, there exists an infinite family of $k$-crossing critical graphs, all having average degree $r$ and crossing number $k$.**

The statement holds for $N_I = $ Ceiling[Max[$f[r_1]$, $f[r_2]$]], where
$f(x) = 240 + \frac{512}{(6-x)^2} + \frac{224}{6-x} + \frac{25}{16\,(x-3)^2} + \frac{40}{x-3}$.